\documentclass[11pt]{article}

\usepackage{amsfonts}
\usepackage{amscd}
\usepackage{amssymb}
\usepackage{amsthm}
\usepackage{amsmath}

 \theoremstyle{plain}
\newtheorem{thm}{Theorem}[section]
\newtheorem{lemma}[thm]{Lemma}
\newtheorem{prop}[thm]{Proposition}
\newtheorem{cor}[thm]{Corollary}
\newtheorem{conj}[thm]{Conjecture}
\theoremstyle{definition}
\newtheorem*{defn}{Definition}

\theoremstyle{remark}
\newtheorem*{remark}{Remark}
\newtheorem*{ack}{Acknowledgments}
\numberwithin{equation}{section}

\def\cA{\mathcal{A}}

\def\cR{\mathcal{R}}

\def\CC{\mathbb{C}}

\def\ZZ{\mathbb{Z}}

\def\fg{\mathfrak{g}}

\def\fsl{\mathfrak{sl}}  
\def\fgl{\mathfrak{gl}}

\def\id{\mathrm{id}}

\makeatletter
\renewcommand{\@makefnmark}{\mbox{\textsuperscript{}}}
\makeatother

\title{Braided central elements}
\author{Natasha Rozhkovskaya}

\date{ }

\begin{document}

\maketitle

\begin{abstract}
We present and study two families of polynomials with coefficients in the center of the universal
enveloping algebra. These polynomials are analogues of a determinant and a characteristic polynomial of a
 certain non-commutative matrix, labeled by irreducible representations of $\fgl_n(\CC)$. The matrix is an image of the universal R-matrix of a Yangian of 
 $\fgl_n(\CC)$ under certain representation. We compute the polynomials  explicitly for $\fgl_2(\CC)$ and establish  connections 
 between the first family of polynomials and higher Capelli identities through some sort of plethysm.

\end{abstract}


\section {Introduction}
Generalizations  of Capelli idenities were studied  in the series of works  by M.\,Nazarov, A.\,Okounkov, A.\,Molev,
 G.\,Olshanskii, F.\,Knop, S.\, Sahi et al (\cite{Naz}, \cite{Mol-Naz}, \cite{Okou1}, \cite{Sahi}, \cite{Knop}, etc.).
 The identities are connected to shifted  symmetric polynomials, satisfy 
some vanishing conditions  and produce a linear basis of the center of the universal enveloping algebra of a Lie algebra.
In \cite{Mol-Naz}, \cite{Naz} these  Capelli identities were studied through the theory of Yanginas.

In this paper we introduce two families of central polynomials, which we believe to have similar links to 
the theory of  representations of Yangian of $\fgl_n(\CC)$. Both series of polynomials are parameterized  by dominant weights
of $\fgl_n(\CC)$. 

One of them can be interpreted as a family of determinants of some non-commutative matrix. We conjecture that
these polynomials have coefficients in the center of the universal enveloping algebra and we prove this 
for  representations of $\fgl_2(\CC)$.
In Section \ref{Rm} we connect "determinants" to Capelli identities through some sort of plethysm.

The second family of central polynomials represents  analogues of characteristic polynomials
 of certain non-commutative matrices. Their existence and centrality
follows from  the works of B.\,Kostant \cite{Kost1} and M.\,Gould \cite{Gould1}.  
We do not know yet the interpretation of these polynomials in terms of Yangians and relations to Capelli elements.  
In case of  the vector representation, "determinant" and "characteristic polynomial" coinside -- both can 
 be obtained as an image of quantum determinant  of the Yangian  of $\fgl_n(\CC)$ under the evaluation map.
 In Section \ref{Sgl2}
we study the case of $\fgl_2(\CC)$, which illustrates that in general these polynomials are different.

  Both series of central polynomials are produced by  certain matrices with coefficients in the universal 
  enveloping algebra. We call these matrices braided Casimir elements. These elements appear in different
  areas of representation theory. For example, in \cite{Kost1} braided Casimir elements were used  to 
  study tensor products of finite and infinite dimensional representations. In \cite{Rozh} it is proved that
   quantum family algebras, introduced by A.\,Kirillov in \cite{Kir}, are commutative if and only if they are
    generated by braided Casimir elements. In \cite{Gould2}, \cite{Gould3} characteristic polynomials of braided 
	Casimir elements were applied to calculate Wigner coefficients. In Section \ref{Rm} 
	 we show that braided Casimir elements are
	images of the universal R-matrix of the Yangian of $\fsl_n(\CC)$ under certain represenations.

\begin{ack} It is my pleasure  to express my gratitude to A.\,Molev and A.\,Ram for
interesting and  valuable discussions.
I would like to thank also   IHES and  IHP for their hospitality and  L.\,V.\,Zharkova
for special support.     	  
\end{ack}

\section{Definitions}\label{Sdef}
 Consider  the universal enveloping algebra 
  $U(\fgl_n(\CC))$ of the  general linear Lie algebra $\fgl_n(\CC)$. Let  $Z(\fgl_n(\CC))$ 
   be the center
  of $U(\fgl_n(\CC))$. 
  Fix the basis $\{E_{ij}\}$ of $\fgl_n(\CC)$, which consists of  standard unit matrices.
   We write  the Casimir element of $U(\fgl_n(\CC))\otimes U(\fgl_n(\CC))$ as
  $$
  \Omega= \sum_{i,j=1,\,\dots,\, n } E_{ij}\otimes E_{ji}.
  $$
  The element $\Omega$ in the tensor square of $U(\fgl_n(\CC))$ is closely related to the central 
  element  $t\in Z(\fgl_n(\CC))$ defined by $t= \sum E_{ij} E_{ji}$
  (here we mean multiplication in the universal enveloping algebra). Namely,
 let $\delta$ be the standard coproduct on  $U(\fgl_n(\CC))$, defined on the 
  elements of $\fgl_n(\CC)$ as $\delta(x)=x\otimes 1+1\otimes x$.
  Then 
  $$
  \Omega=\frac{1}{2}\left(\delta(t)-1\otimes t -t\otimes 1
  \right). $$

Let   $\lambda=\, (\lambda_1,  \dots , \, \lambda_n)$ with $\lambda_i-\lambda_{i+1} \in \ZZ_+$, for 
$i=1,\,\dots,\,n-1$, be  a dominant weight of $\fgl_n(\CC)$.
 Denote by  $\pi_\lambda$ be the corresponding 
irreducible  rational $\fgl_n(\CC)$-representation  and by $V_\lambda$
 the space of this representation. We assume that 
 dim $V_\lambda=m+1$. 	
  We construct 
from $\Omega$ a new object, an element of
   $ U(\fgl_n(\CC))\,\otimes \,\text{End}(V_\lambda)$, which  we  call {\it braided Casimir element}.

  \begin{defn} The braided Casimir element is  defined by
  $$
  \Omega_\lambda=E_{ij}\otimes\sum_{i,j=1,\,\dots ,\,n} \pi_\lambda(E_{ji}).
  $$
  \end{defn}
  We will  think of $\Omega_\lambda$ as a matrix of size $(m+1)\times (m+1)$ 
  with   coefficients in  $U(\fgl_n(\CC))$.  
   Hence,  braided Casimir element comes from the central element $t$ of the universal enveloping algebra. It turns out that
   $\Omega_\lambda$  is a source of 
new central  elements.  Below we define   two families  of  polynomials with coefficients in $Z(\fgl_{n}(\CC))$, associated 
to $\Omega_\lambda$.
 
 The first family of polynomials
	 is provided by the following proposition.
	 
	 Let $\fg$ be a simple Lie algebra with a  linear basis $\{I_\alpha\}$. Let $\{I^\alpha\}$ be the 
	 dual basis with respect to Killing form.
Consider an element $\omega$ of $U(\fg)\otimes U(\fg)$, defined similarly by $\omega=\sum_{\alpha} I_\alpha\otimes I^\alpha$.
For any irreducible representation $\pi$ of $\fg$ put
$\omega_\pi=(\id\otimes \pi)(\omega) $. 
	 
\begin{prop}\label{ss}
(\cite {Kost1},\cite{Gould1})
There exists a polynomial $$p_\pi(u)=\sum_{k=0}^{m}z_{k}u^k$$  with coefficients
$z_k$ in the center of $U(\fg)$ such that
 $
p_\pi(\omega_\pi)=0$.

The polynomial $p_\pi(u)$  can be chosen so that  $\text{deg}\, p_\pi(u)=\, \text{dim}\, \pi$. 
\end{prop}	   
\begin{cor} For any  dominant weight $\lambda$ of $\fgl_n(\CC)$ 
there exists a polynomial 
$$P_\lambda(u)=\sum_{k=0}^{m}z_k u^k$$  with coefficients 
$z_k\in Z(\fgl_n(\CC))$, such that $P_\lambda(\Omega_\lambda)=0$.
\end{cor}
\begin{proof} Choose a basis in $\fsl_n(\CC)$ which consists of matrices
$E_{ij}$ for $i\ne j$ and $H_i=E_{ii}-E_{i+1,i+1}$,  $(i=1,..., n-1)$.
The dual basis consists of matrices 
$$E^*_{ij}=E_{ji}\quad \quad {\text and} 
\quad \quad H^*_i=(E_{11}+...+E_{ii})-\frac{i}{n}(E_{11}+...+E_{nn}).$$ 
The element  $\omega$ for $\fsl_n(\CC)$  has the form
$$
\omega=\sum_{i\ne j} E_{ij}\otimes E_{ji}+\sum_i{H_i}\otimes H^*_i=
\sum_{ij}E_{ij}\otimes E_{ji} - \frac{1}{n}(\Delta_1\otimes \Delta_1)=\Omega-\frac{1}{n}
(\Delta_1\otimes \Delta_1),
$$
where $\Delta_1=(E_{11}+ \dots +E_{nn})$ is a central element of $U(\fgl_n(\CC))$.
Consider a dominant
 weight $\lambda^\prime$ of $\fsl_n(\CC)$, defined as $\lambda^\prime=(\lambda_1-\lambda_m,\,\dots , \,\lambda_{m-1}-\lambda_m, \,0)$. 
 Then
 $$
 \omega_{\pi_{\lambda^\prime}}=\Omega_\lambda - \frac{d}{n}\Delta_1\otimes id
 $$ 
 with $d=\sum_{i}\lambda_i$,
and the polynomial $p_{\lambda^\prime}(u- \frac{d}{n}\Delta_1)$ from Proposition \ref{ss}
annihilates $\Omega_\lambda$. 
 \end{proof}
 We will call  $P_\lambda(u)$ {\it characteristic polynomial} of $\Omega_\lambda$.
   
   Next we define the second family of polynomials $D_\lambda(u)$, which  we  call
   {\it shifted determinants}.
   
  Let $A$ be an element of  
     $ \mathcal A\otimes\,\text{End}\,(\CC^{m+1})$, where 
	 $\cA$ is a non-commutative algebra, let 
   $V=\CC^{m+1}$ be an $(m+1)$-dimensional vector space.
	  We again think of $A$ as
	a  non-commutative matrix of size  $(m+1)\times (m+1)$ with coefficients
	$A_{ij}\in \mathcal A$. 
  \begin{defn} 
  The (column)-determinant of  $A$ is the following element of $\mathcal A$:
  \begin{equation}
  \text{det}(A)=\sum_{\sigma\in S_{m+1}}(-1)^{\sigma}A_{ \sigma(1)1}A_{\sigma(2)2}... A_{ \sigma(m+1)(m+1)}.
  \end{equation}
  Here the sum is taken over all elements $\sigma$ of the symmetric group $S_{m+1}$ and $(-1)^\sigma$
   is the sign of the permutation $\sigma$.
  \end{defn}
  Put $\Omega_\lambda(u)=\Omega_\lambda+u\otimes id$.
   Define  $L$ as a  diagonal matrix
   of the size $(m+1)\times (m+1) $ of the  form:  
   $$
   L=\text{diag}(m,m-1,\, \dots, \,0).
   $$
\begin{defn}
 The shifted determinant of $\Omega_\lambda(u)$ is  
 the column-determinant $\text{det} (\Omega_\lambda(u)-L)$.
 We will use notation $D_\lambda(u)$ for this polynomial with coefficients in $U(\fgl_n(\CC))$:
 $$D_\lambda(u)=\text{det} (\Omega_\lambda(u)-L).$$ 
 \end{defn}  	   
\begin{conj} \label{conj}
For any dominant weight $\lambda$  there exists a basis of the vector space $V_\lambda$ such that  the polynomial
 $D_\lambda (u)$ has coefficients in the center $Z(\fgl_n(\CC))$.
\end{conj}  
 This is known to be true in the case of vector representation and it is proved below in Section
 \ref{Sgl2} for  representations of $\fgl_2(\CC)$. 
\smallskip

There is another way to define  the same determinant.  	   
Let $A_1,\,\dots,\, A_s$ be a set of matrices of size $(m+1)\times (m+1)$  
with coefficients in some associative (non-commutative) algebra $\mathcal A$.
 Let $\mu$ be  the multiplication in $\mathcal A$.
Consider an element of $\mathcal A\otimes \text{End}(V)$
$$
\Lambda^{s}(A_1\otimes \dots \otimes A_s)= 
(\mu ^{(s)}\otimes Asym_{s})(A_1\otimes\dots\otimes A_s),
$$
where $Asym_{s} =\frac{1}{s!}\sum_{\sigma\in S_s }(-1)^{\sigma}\sigma$. By Young's construction, the antisymmetrizer
 can be  realized as an element of End $(V^{\otimes s})$.
\begin{lemma}(cf \cite{MNO}.) For $s=m+1= \text{dim}\, V $
\begin{equation}\label{asym}
\Lambda^{m+1}(A_1\otimes \dots\otimes A_{m+1})=\alpha(A_1,\dots, A_{m+1})\otimes Asym_{m+1}, 
\end{equation}
where
$\alpha(A_1,\dots, A_{m+1})\in \mathcal A$,
$$
\alpha(A_1,\dots, A_{m+1})=
\sum_{\sigma\in S_{m+1} }(-1)^{\sigma}[A_1]_{\sigma(1),1}\dots [A_{m+1}]_{\sigma (m+1),m+1},
$$
$[A_k]_{i,j}$ -- matrix elements of $A_k$.
\end{lemma}
\begin{proof}Let $\{e_i\}$,  $(i=1,\,\dots,\, m+1)$,  be a basis of $V$.
Observe that $Asym_{\,(m+1)}$ is a one-dimensional projector to 
$$
v=\frac{1}{(m+1)!}\sum_{\sigma\in S_{m+1}}(-1)^{\sigma}e_{\sigma(1)}\otimes\dots \otimes e_{\sigma(m+1)}.
$$
We  apply  $\Lambda^{m}(A_1\otimes \dots \otimes A_{\,m+1})$ to $e_1\otimes\dots\otimes e_{m+1} \in V^{\otimes m}$: 
\begin{equation}\label{a}
\begin{split}
\Lambda^{m}(A_1\otimes \dots \otimes A_{\,m+1})(e_1\otimes\dots \otimes e_m)& \quad \\
=\sum_{i_1,\dots i_k} (A_1)_{\,i_1,\,1}\dots (A_{\,m+1})_{\,i_{m+1},\,m+1} 
\,Asym_{\,(m+1)}\, 
(e_{i_1}\otimes \dots \otimes e_{i_{m+1}}). &\\
\end{split}
\end{equation}
The vector $Asym_{\,(m+1)}\,(e_{i_1}\otimes \dots \otimes e_{i_{m+1}})\ne 0$ only if all indices $\{i_1,\,\dots,\, i_{m+1}\}$
 are pairwise distinct.
In this case denote by $\sigma$ be a permutation defined by $\sigma(k)=i_k$. Then  
 $$Asym_{\,(m+1)}\,(e_{i_1}\otimes \dots \otimes e_{i_{m+1}})=(-1)^{\sigma}\,v,$$
  and (\ref{a}) gives 
\begin{equation*}
\begin{split}
\Lambda^{m}(A_1\otimes\dots \otimes A_{m+1})\,(e_1\otimes\dots \otimes e_{m+1})=
\alpha(A_1,\,\dots,\, A_{\,m+1})\,v\\
=\alpha(A_1,\,\dots,\, A_{\,m+1})\,Asym_{\,(m+1)}\,(e_1\otimes\dots \otimes e_{m+1}).
\end{split}
\end{equation*}
\end{proof}
\begin{remark}
Note that
$$ 
\alpha \,(A,\,\dots,\, A)=\text{det}(A),
$$
$$
\alpha\, (\Omega_\lambda(u-m),\, \dots, \, \Omega_\lambda(u))=D_\lambda(u).
$$
\end{remark}

\section{Yangian of $\fgl_n (\CC)$  and  Casimir element}\label{yng}
In this section we would like to recall  the case of vector representation and its connection to  Yangians.
This example   serves  as an inspiration for  the rest of the project. 
Let us recall some definitions (\cite{Dr1}, \cite{MNO}).
\begin{defn} The Yangian $Y(n)$ for $\fgl_n(\CC)$ is a unital  associative algebra over $\CC$ with 
countably many generators $\{t_{ij}^{(r)}\}$, $r=1,2, ...$, $1\le i,j\le n$ and the defining 
relations 
$$
[t_{ij}^{(r+1)},t_{kl}^{(s)}]-[t_{ij}^{(r)},t_{kl}^{(s+1)}]=t_{kj}^{(r)}t_{il}^{(s)}-t_{kj}^{(s)}t_{il}^{(r)},
$$
 where $r,s=0,1,2...$ and $t_{ij}^{(0)}=\delta_{ij}$.
\end{defn}
The same  set of defining relations can be combined into  one equation, sometimes called {\it RTT relation}.
Namely, denote by  
$T(u)=(t_{ij}(u))_{i,j\, =1}^{n}$  the matrix with coefficients $t_{ij}(u)$ which a formal power  series of
 generators of  $Y(n)$ : 
$$t_{ij}(u)=\delta_{ij}+\sum_{k=1}^{\infty} \frac{t_{ij}^{(k)}}{u^k}. $$
For $
P={\sum E_{ij}\otimes E_{ji}}
$,  the permutation matrix of $\CC^n\otimes \CC^n$, define the Yang matrix
$$R(u)=1-\frac{P}{u}.$$

$R(u)$ is a rational function with values in $\text{End}\, \CC^n\otimes  \text{End}\, \CC^n $.

We introduce some standard notations. 
For any vector space $V$ and any  element $S$ of $\text{End}\, V $
 we define an element $S_k$ of $\text{End} V^{\otimes m} $
by 
$$
S_{k}=1^{\otimes(k-1)} \otimes S\otimes 1^{\otimes(m-k)}.
$$
In particular, we write
$$
T_{k}(u)=\sum_{ij} t_{ij}(u)\otimes (E_{ij})_{k} \, \in Y(n)\otimes\text{End} (\CC^n)^{\otimes m}.
$$
  Let $S$ be an element of $\text{End} V\otimes \text{End} V $.
  Using  the abbreviated notation $S=S(1)\otimes S(2)$,
    we define an element $S_{ij}$ of $\text{End} (\CC^n)^{\otimes m} $
by 
$$
S_{ij}=1^{\otimes(i-1)}\otimes S{(1)}\otimes 1^{\otimes (j-i-1)}\otimes S{(2)} \otimes 1^{\otimes (m-j-i)}.
$$
 
\begin{defn} The Yangian $Y(n)$ of $\fgl_n(\CC)$ is an associative unital  algebra over $\CC$ with the set of generators
$\{t_{ij}^{(k)}\}$ which satisfy the equation
$$
R(u-v)T_1(u)T_2(v)=T_2(v)T_1(u)R(u-v).
$$
 (This is an equation in $Y(n)\otimes \text{End}\, (\CC^n)^{\otimes 2}$). 
\end{defn}
Yangian $Y(n)$  is an example of infinite-dimensional quantum group.
 It has a remarkable  group-like central element, which is called {\it quantum
 determinant} of $Y(n)$.
\begin{defn} Quantum determinant $\text{qdet}\,T(u)$ is a formal series with coefficients in $Y(n)$
defined by
$$\text{qdet} T(u)=\sum_{\sigma\in S_n} (-1)^\sigma t_{1\sigma(1)}(u-n+1) \dots t_{n\sigma(n)}(u).$$
\end{defn}

  
  Let $V$ be a $\fgl_n(\CC)$-module. The action $\pi$ of $\fgl_n(\CC)$ on $V$ can be extended to the action of the Yangian
   $Y(n)$. Usually it is done by the means of evaluation map which is a homomorphism of algebras
   $\text{ev}:\,Y(n)\, \to\, U(\fgl_n(\CC))$. By definition,
  
 $$ \text{ev}\cdot t_{ij}(u)=\delta_i^j+\frac{E_{ij}}{u}.$$
 It is clear that 
 $$
\text{ev}\cdot T(u)=1+\frac{\Omega^\top}{u},
$$
where $\top$ is for  the matrix transposition.
Along with $(\text{ev})$ we will consider another homomorphism $\check{\text{ev}}:\, Y(n)\,\to\, U(\fgl_n(\CC))$, 
which is defined by
 $$\check{ \text{ev}}\cdot t_{ij}(u)=\delta_i^j-\frac{E_{ji}}{u},\quad \quad
\check{\text{ev}}\cdot T(u)=1-\frac{\Omega}{u}.$$
In Section \ref{Rm} we discuss these two maps in more details.

One can see that   
$$
\text{ev}(\text{qdet}T(u))=\frac{ D_{\lambda_0}(u)}{u(u-1)...(u-n+1)},
$$
where $D_{\lambda_0}(u)$ is  the "shifted"  determinant for the  vector representation $\lambda_0=(1,0,\,...,\, 0)$ :
\begin{equation*}
D_{\lambda_0}(u)=\text{det}
\begin{pmatrix}
 E_{11}+u-n+1 & E_{21}& ...& E_{n,1}\\
E_{12}& E_{22}+u-n+2&...& E_{n,2}\\
...& ...& ...& ...\\
 E_{1n}& E_{2n}& ...& E_{n,\,n}+u
\end{pmatrix}.
\end{equation*}
\begin{prop}(\cite{Naz-Tar} \cite{Mol1}.)
$D_{\lambda_0}(u)$ is a polynomial with coefficients in the center of $U(\fgl(n))$. Moreover,
for vector representation  $$D_{\lambda_0}(-u)=P_{\lambda_0}(u),$$ -- the shifted  determinant
 provides the characteristic polynomial
of the corresponding braided Casimir element $\Omega_{\lambda_0}$.
\end{prop}
As we will see in Section \ref{Sgl2}, these   two polynomials do not coincide in general.

  \section{Case of $\fgl_2(\CC)$}\label{Sgl2}
  Here we prove the centrality of polynomials $D_\lambda(u)$ 
  and compare  $D_\lambda(u)$ and $P_\lambda(u)$ for irreducible representations of $\fgl_2(\CC)$. 
  The center $Z(\fgl_2(\CC))$ of the universal enveloping algebra $U(\fgl_2(\CC))$ is
   generated by two elements:
   $$\Delta_1=E_{11}+E_{22}, \quad \quad \Delta_2=(E_{11}-1)E_{22}-E_{12}E_{21}.$$
  Let $\lambda=(\lambda_1 \ge \lambda_2) $ be a dominant weight. Put $m=\lambda_1-\lambda_2$,\quad $d=\lambda_1+\lambda_2$.
  Then $\text{dim}V_\lambda=m+1$ and 
  $\Omega_\lambda$ is a "tridiagonal" matrix: all  entries $[\Omega_\lambda]_{ij}$ of the matrix $\Omega_\lambda$
   are zeros except
  \begin{align*}
  [\Omega_\lambda]_{k,k}=(\lambda_1-k+1)E_{11}+(\lambda_2+k-1) E_{22},& \quad 
  \quad k=1,\, \dots,\, m+1, \\
  [\Omega_\lambda]_{k,k+1}=(m+1-k)E_{21},&  \quad\quad k=1,\,\dots, \,m,&\quad\\
  [\Omega_\lambda]_{k+1,k}=kE_{12},&\quad  \quad k=1,\,\dots, \,m.&\quad
  \end{align*}

  \begin{prop} a) Polynomial $D_\lambda (u)$ is central.
  
  b) Let $\mu=(\mu_1\,\ge\, \mu_2 )$ be another dominant weight of $\fgl_n(\CC)$.
  The image of $D_\lambda (u)$ under Harish - Chandra isomorphism $\chi$  
  is the following function of $\mu$:
  \begin{equation}\label{D-HC}
  \chi(D_\lambda(u))=\prod_{k=0}^{m}(u+(\lambda_1-k)\mu_1+(\lambda_2+k)\mu_2 -k).
  \end{equation}
  \end{prop}
  \begin{proof} a) 
 Let $X=X(a,b,c)$ be a matrix of size $(m+1)\times (m+1)$
of the form
$$
\begin{pmatrix}
a_m&b_m&0&\dots&0&0&0\\
c_m&a_{m-1}&b_{m-1}&\dots&0&0&0\\
\dots&\dots&\dots&\dots&\dots&\dots&\dots\\
0&0&0&\dots&c_2&a_{1}&b_{1}\\
0&0&0&\dots&0&c_1&a_{0}\\
\end{pmatrix},
$$
where $a_{i}, b_{j}, c_{k}$ are elements of some (noncommutative)
algebra. 
 Define $\text{det}\,X$ as in Section
\ref{Sdef}. Using
the principle $k$-minors
of $X$, the determinant of $X$ can be computed by recursion formula.
Denote 
by  $X_k$ the matrix obtained from
$X$ by deleting the first $(m-k+1)$
rows and the first $(m-k+1)$ columns, and by $I^{(k)}$  the determinant of $X_k$.
 Then det\,$X=I^{(m+1)}$, and
we have the following recursion:
\begin{lemma} The determinants $I^{(k)}$ satisfy the recursion  realtion
\begin{equation}
\begin{split}
I^{(k+1)}=a_kI^{(k)}-c_kb_k\,I^{(k-1)},\quad \quad (k=2,...,n-1)
\end{split}
\end{equation}
{with initial conditions}
$I^{(1)}=a_0$, 
$I^{(2)}=a_1a_0-c_1b_1$.
\end{lemma}
\begin{lemma}\label{subs}
If $ X(a,b,c)$ is a tri-diagonal matrix with coefficients
$\{a_i,b_j,c_k\}$
and $ X(a',b',c')$ is another tri-diagonal matrix with
coefficients $\{a'_i,b'_j,c'_k\}$
with the property
$$
a_i=a^\prime _i, \quad \quad c_j\,b_j=c^\prime _j\, b^\prime _j
$$
for all $i=0,...,n$ $j=1,...,n$,
then
$ \text{det}\,X(a,b,c)=\text{det}\, X(a^\prime, b^\prime,
c^\prime)$. 
\end{lemma}
\begin{proof}
Follows from the recursion relation.
\end{proof}
 We apply these observations   to compute the determinant of the matrix
$X(a,b,c)=\Omega_\lambda(u)-L$
with parameters
$$a_k=(\lambda_1+k-m)E_{11}+(\lambda_2+m-k)E_{22}+u-k, \quad\quad  k=0,\,\dots,\,m, $$
$$c_k=(m-k+1)\, E_{12}, \quad\quad\quad  b_k=kE_{21},\quad \quad k=1,\,\dots,\,m.$$

The following obvious lemma allows to reduce the  determinant of non-commutative matrix $(\Omega_\lambda(u)-L)$
to a  determinant of a matrix with commutative coefficients.
\begin{lemma}\label{com}
 The subalgebra of $U(\fgl_2(\CC))$ generated by
$\{E_{11}, E_{22}, (E_{12}\,E_{21}) \}$ is commutative.
\end{lemma}
Put
$h=E_{11}-E_{22}$, 
 $a=
E_{12}E_{21}$.
Due to Lemma \ref{subs}  the tridiagonal  matrix  $X(a^\prime,b^{\prime},c^{\prime})$
with coefficients
\begin{equation*}\begin{split}
a^\prime_k=a_k=\lambda_1E_{11}+\lambda_2E_{22}+u-m+(k-m)(h-1),& \quad \quad k=0,\,\dots,\,m,\\
b^\prime_k=k\,a, \quad \quad \quad
c^\prime_k=(m-k+1),&\quad \quad k=1,\,\dots,\,m,
\end{split} 
\end{equation*}
 has the same determinant as $(\Omega_\lambda(u)-L)$. By Lemma \ref{com}, 
$X(a^\prime,b^{\prime},c^{\prime})$ has commutative coefficients.
Hence 
$\text{det}\,(\Omega_\lambda(u)-L)$ equals 
$\text{det}\,(\lambda_1E_{11}+\lambda_2E_{22}+u-m+ A_m)$
where $A_m$ is the following matrix:
$$A_m=\begin{pmatrix}
0&ma&0&\dots&0&0\\
1&-(h-1)&(m-1)a&\dots&0&0\\
\dots&\dots&\dots&\dots&\dots&\dots\\
0&0&0&\dots&(1-m)(h-1)&a\\
0&0&0&\dots&m&-m(h-1)\\
\end{pmatrix}
$$
with  $h$ and $a$  as above.

\begin{lemma}
 $$\text{det}\, A_m
=\prod_{k=0}^{m}\left(\frac{-m(h-1)}{2}+
\frac{(m-2k)}{2}{\sqrt{(h-1)^2+4a}}\right).
$$
\end{lemma}
\begin{proof} By Lemma \ref{subs}
$\text{det}\, A_m={(h-1)^{m+1}} \text{det}\, A^\prime_m$
with
$$
A^\prime_m=
\begin{pmatrix}
0&ms&0&\dots&0&0\\
1(s-1)&-1&(m-1)s&\dots&0&0\\
\dots&\dots&\dots&\dots&\dots&\dots\\
0&0&0&\dots&1-m&s\\
0&0&0&\dots&m(s-1)&-m\\
\end{pmatrix},
$$
and $s$ is such that $s(s-1)={a}/{(h-1)^2}$. This reduces to
\begin{equation} \label {s}
s=\frac{1\pm \sqrt{1+4a/(h-1)^2}}{2}.
\end{equation}
The determinant of $A^\prime_m$ is a variant of Sylvester
determinant (\cite{Askey}, \cite{Holtz}, \cite{Sylv}).
It equals
$$\text{det}\,A_m^\prime
=\prod_{k=0}^{m}((m-2k)s-m+k).
$$
With $s$ as in (\ref{s})
we have:
\begin{equation}\notag
\begin{split}
(h-1)((m-2k)s-m+k)
=\frac{-m(h-1)}{2}\pm\frac{(m-2k)}{2}\sqrt{\left({(h-1)^2+4a}\right)}
\end{split}
\end{equation}
and lemma follows. Note that both values of $s$ give the same value of  det $A_m$.
\end{proof}
We obtain from calculations above
\begin{equation} 
\text{det}\, (\Omega_\lambda(u)-L)=\prod_{k=0}^{m}\left(u+\frac{d}{2}(E_{11}+E_{22})
-\frac{m}{2}+\frac{(m-2k)}{2}{\left({(E_{11}-E_{22}-1)^2+4E_{12}E_{21}}\right)^{\frac{1}{2}}}\right).
\end{equation}
Observe that ${\left({(E_{11}-E_{22}-1)^2+4E_{12}E_{21}}\right)^{\frac{1}{2}}}=\left({
(\Delta_1-1)^2-4\Delta_2}\right)^{\frac{1}{2}}$,
and finally we get
\begin{equation}\label {O_n}
D_\lambda(u)=\prod_{k=0}^{m}
\left(u+\frac{d\,\Delta_1}{2}
-\frac{m}{2}+\frac{(m-2k)}{2}\left({(\Delta_1-1)^2-4\Delta_2}\right)^{\frac{1}{2}}\right).
\end{equation}
 The quantity  in
(\ref{O_n}) has coefficients in $W_\tau$-extension of $Z(\fgl_2(\CC))$, where $W_\tau$
is the translated Weyl group. 
But it is easy to see that after expanding the product, we get a polynomial in $u$ with coefficients
in $Z(\fgl_n(\CC)$. We proved the first part of the proposition.

b) The images of the generators of $Z(\fgl_n(\CC))$ under Harish-Chandra  homomorphism $\chi$ are
$$
\chi(\Delta_1)=\mu_1+\mu_2,\quad \quad\chi(\Delta_2)=\mu_1(\mu_2-1).
$$   
This together with (\ref{O_n}) implies (\ref{D-HC}).
 \end{proof}
 
 \begin{prop}
 The image of the polynomial $P_\lambda(u)$ under Harish-Chandra homomorphism is 
 \begin{equation}\label{P-HC}
  \chi(P_\lambda(u))=\prod_{k=0}^{m}\left(-u+(\lambda_1-k)\mu_1+(\lambda_2+k)\mu_2 -k(m+1-k)\right).
  \end{equation}
 \end{prop}
 This follows from the formula for characteristic polynomial for  $\omega_n$, the braided Casimir element of $\fsl_2(\CC)$.
 It is proved in \cite{Rozh}.
   
   Despite the fact that in general $D_\lambda(u)$ and $P_\lambda(u)$ are different, we believe that
 still they are closely related. We hope to find  the source of the link between these two polynomials in 
 representation theory of Yangian $Y(n)$.
  Some steps towards this goal are presented in the next section.

    \section{Connection with Yangian $Y(\fgl_N(\CC))$}\label{Rm}

 As  we saw in the Section \ref{yng},  the Casimir
 element $\Omega$ has a natural interpretation as an image of the matrix of generators of the Yangian
  under the  map $\check{\text{ev}}$.
  Here we give some similar interpretation of the  braided Casimir element. Though, in this interpretation
  the element $\Omega_\lambda$ arises not as an element of $U(\fgl_n(\CC))\otimes \text{End}\, V_\lambda$, but
  as an element of $\text{End}\,\CC^n \otimes \text{End}\, V_\lambda$.
  
  Let  $\cR(u)$ be the universal R-matrix of the Yangian  of $\fsl_n(\CC)$  (see \cite{Dr1}). It is a formal power series
  in $u^{-1}$ with coefficients in $Y(n)\otimes Y(n)$.
  
  Let $\pi$ be a $\fgl_n(\CC)$-representation in  some space $V$. Using the maps
   $(\text{ev})$ and $(\check{\text{ev}})$
  we can construct two actions $\pi$ and $\check{\pi}$ of $Y(n)$ on $V$:
  $$
  \pi\,\cdot\, y =(\pi\circ \text{ev})\,\cdot\, y  \quad \quad 
  \text{and} \quad 
    \check{\pi}\,\cdot\, y =(\pi\circ \check{\text{ev}})\,\cdot \, y
  $$  
  We need both actions to make the correspondence between the images of the universal R-matrix,
  braided Casimir element, and Capelli polynomials.
\begin{remark}
It is well-known that
the universal R-matrix satisfies the quantum Yang-Baxter equation with spectral parameter:
\begin{equation}\label{QYBE}
\cR^{12}(u-v)\cR^{13}(u-w)\cR^{23}(v-w)=\cR^{23}(v-w)\cR^{13}(u-w)\cR^{12}(u-v).
\end{equation}
Let $\tilde\pi$ be some representation of $Y(n)$, extended in some way from the representation $\pi$ of $U(\fgl_n(\CC))$. After the application of $\tilde\pi\otimes\tilde\pi\otimes id$ 
to both sides of (\ref{QYBE}) we get
exactly the RTT-relation for the matrix $T(u)=(\tilde \pi\otimes id)\cR(u)$ of elements in $Y(n)$
with the defining R-matrix $(\tilde\pi\otimes \tilde\pi)\cR(u)$.
We would like to stay with the tradition and we would like to  get in the case of 
$\tilde\pi$, extended from the
 vector representation of $\fgl_n(\CC)$,
the Young R-matrix $1-P/u$ and the  defining relation of $Y(n)$.
Observe that {for any } $\tilde\pi$
$$
(\id\otimes \tilde\pi) T(u)= (\tilde\pi\otimes \tilde\pi) \cR(u).
$$  
 But  for the  vector representation $\pi_0$ of $\fgl_n(\CC)$  and  the matrix $T(u)$ of generators of 
 $Y(n)$ we have 
 $$ (\id\otimes \pi_0) T(u)= 1+\frac{P^\top}{u},\quad \quad 
 (\id\otimes \check{\pi}_0)T(u)=1-\frac{P}{u}=R(u).
 $$
Hence, to be consistent in definitions, we have to extend the representations from 
$U(\fgl_n(\CC))$ to $Y(n)$ by the map $(\check{\text{ev}})$:
\begin{equation}\label{R_v}\tilde\pi=\check\pi,\quad\quad
R(u)=f_0(u)(\check{\pi_0}\otimes\check{\pi_0})\cR(u), \quad \quad  
T(u)=(\check{\pi_0}\otimes id)\cR(u),
\end{equation}
for some complex-valued  rational function $f_0(u)$.
 \end{remark}
  \begin{prop} \label{R_im}
  Let $\lambda_0=(1,\,0, ...,\, 0)$ be the highest weight of the vector representation $\pi_0$
  of $\fgl_n(\CC)$, let $\lambda=(\lambda_1\,\ge\, \lambda_2\, \ge\,\dots\,\ge \lambda_n\, \ge\, 0)$ be a partition (a  dominant weight of $\fsl_n(\CC)$).
  Let $\check{\pi}_0$, $\check{\pi}_\lambda$ be the representations of $Y(n)$, extended from $\fgl_n(\CC)$ by 
  the map $(\check{\text{ev}})$. 
     Then$$
(\check{\pi}_{0}\otimes \check{\pi}_\lambda)\cR(u)= f(u)\left(1-\frac{\Omega_\lambda}{u}\right),
$$
where $f(u)$  is some complex-valued  rational function of $u$ (which depends on $\lambda$).
\end{prop}

\begin{proof}
In general, this follows from the Corollary 3.\,7 in \cite{Naz}. We 
repeat the steps which would lead to  the proposition to avoid confusion
 in notations. 
Let $V_0=\CC^n$ be the  space of the  vector representation  $\pi_0$ with the highest weight 
$\lambda_0$. By the means of the  map $\check{ev}$ it becomes $Y(n)$-module.
Moreover, there exists a family of $Y(n)$-modules $\{V_0(a)\}_{a\in \CC}$. The action 
$\check{\pi}_{0}(a)$ of $Y(n)$ on $V_0(a)$ is defined as follows:
$$
\check{\pi}_{0}(a)\cdot T(u)=\check{\pi}_{0}  \cdot T(u-a),
$$ 
where $T(u)$ is the matrix of generators of $Y(n)$.

 For two dominant weights $\nu$ and  $\lambda$ of $\fsl_n(\CC)$ put 
$\cR_{\nu,\,\lambda}(u)=(\check{\pi}_{\nu} \otimes \check{\pi}_\lambda)\cR(u)$. 
By the discussion  (\ref{R_v}) above, the image of the universal R-matrix under the representation
$\check{\pi}_{0}\otimes \check{\pi}_{0}$  is proportional to the Yang R-matrix.
It is also well-known that the universal R-matrix satisfies
\begin{equation}\label{delta}
(id\otimes\Delta)(\cR(u))=\cR_{12}(u)\cR_{13}(u),
\end{equation}
where $\Delta$ is a coproduct in $Y(n)$.
Let $\lambda$ be a partition of $M$ and let $\{ c_1,\,\dots,\, c_n \}$ be the set of 
contents of the standard  Young tableau of the shape $\lambda$. 
Using (\ref{delta}) we can expand the action of the universal R-matrix  to the 
$Y(n)-$module $W=V_0(0)\otimes V_0(c_1)\otimes\dots \otimes V_0(c_M)$,
$i=1 ,\, \dots,\,  M$. This action is 
\begin{equation}\label{tensor}
\check{\pi}_{0}\otimes\check{\pi}_{0} (c_1)\otimes\,\dots\,\otimes \check{\pi}_{0}(c_M)\,
\cR(u)=f(u)\prod_{k=1}^{M}\left(1-\frac{P_{1,k+1}}{u-c_k}\right),
\end{equation}
with $f(u)$ -- some rational complex-valued function of $u$.
 Observe that $W=(V_0)^{\otimes (M+1)}$ as $\fgl_n(\CC)$-module. 
 Using the Young symmetrizer $F_\lambda$ we construct
 the element $(id\,\otimes\, F_\lambda)$ of $\text{End}\,(W)$, which is a projector from $W$ to 
 $V_0\otimes V_\lambda$. We apply this projector to (\ref{tensor}) to obtain 
 $\cR_{\lambda_0\, \lambda}(u)$.
 By Proposition 2.12 in \cite{Naz}, the following equality holds:
 \begin{equation}\label{Naz}
 \prod_{k=1}^{M}\left(1-\frac{P_{1,\,k+1}}{u-c_k}\right) (id\,\otimes\,F_\lambda)
 =\left( 1- \sum_{k=1}^{M}\frac{P_{1,\,k+1}}{u} \right)(id\,\otimes\,F_\lambda).
 \end{equation}     

  But  with the standard coproduct $\delta$ in $U(\fgl_n(\CC))$ we obtain:
\begin{align}
\Omega_\lambda=\sum_{ij}E_{ij}\otimes \pi_\lambda(E_{ji})=
\left(\sum_{ij}E_{ij}\otimes \delta^{(M)}(E_{ji})\right)(id\otimes F_\lambda)
=\left(\sum_{l=1...M,}P_{1,l+1}\right)(id\otimes F_\lambda),\label{omega}
\end{align}
and we get
$$
\cR_{\lambda_0 \, \lambda}(u)=f(u)\left( (id\otimes F_\lambda)- \frac{\Omega_\lambda}{u}\right).
$$
The operator $(id\otimes F_\lambda)- \frac{\Omega_\lambda}{u}$
  acts as  $1-\frac{\Omega_\lambda}{u}$ on the subspace $V_0\otimes V_\lambda \subset W$.
  \end{proof}
 \section{Capelli elements.}

 Next we would like to show how Capelli elements fit into this picture and to 
connect   the  polynomial $D_\lambda(u)$ to these elements  by  some sort of plethysm.
Let  $S=\sum_{ij} E_{ij}\otimes E_{ij}$.
 Following \cite{Okou1},  define an element $S_\lambda(u)$
  of $U(\fgl_n(\CC)\otimes \text{End}V_\lambda$ by 
 $$
 S_\lambda(u)=\left((S_{12}-u-c_1)\, \dots \, (S_{1\,M+1}-u-c_M)\right)\,(id \otimes F_\lambda).
 $$
 Then  $$c_\lambda(u)= tr\,(S_\lambda(u))$$ is {\it the Capelli polynomial}, associated to $\lambda$.
 It has coefficients in $Z(\fgl_n(\CC))$. The theory of Capelli elements is developed in full in the papers, mentioned in the Introduction.
 Here we would like to  observe the following facts:
 
 1)  $S_\lambda(u)$  is proportional to
  $(\pi_0\otimes \check{\pi}_\lambda)\cR(u)$
   (the  first representation  is extended by evaluation map $(\text{ev})$, and the
 second one  by $(\check{\text{ev}})$).
 
 2) Let $\phi:\,\fgl_n(\CC)\to \fgl_n(\CC)$ be an automorphism, defined by $\phi(X)=-X^{\top}$. Put
$\pi_{\lambda^{\star}}=\pi_\lambda\,\circ\,\phi$.  
Then 
\begin{equation}\label{Sl}
 \Omega_{\lambda^\star}(u)= u \prod_{k=1}^{M}\frac{(-1)}{(u+c_k)}\, S_\lambda(u).
\end{equation}
 Indeed, from (\ref{omega}) and (\ref{Naz}) we get
\begin{equation}\label{omega1}
\Omega_\lambda(u)=u\,\prod_{k=1}^{M}\left(1+\frac{P_{1,\,k+1}}{u+c_k} \right)\,(id\otimes F_\lambda).
\end{equation}
 Observe that 
 $(\phi\otimes \pi_\lambda)\,\Omega=id\otimes (\pi_{\lambda}\circ \phi)\, \Omega$,
so
\begin{equation}\label{omega*}
\Omega_{\lambda^\star}(u)=
u\,\prod_{k=1}^{M}\left(1-\frac{S_{1,\,k+1}}{u+c_k} \right)\,(id\otimes F_\lambda).
\end{equation}
and (\ref{Sl}) follows.
 \smallskip 
 
 The representation $X \to - (\pi_{\lambda^\star}(X))^\top$,  $X\in \fgl_n(\CC)$ is isomorphic
 to $\pi_\lambda$ (it has  the same set of weights). 
 Thus we can write in some basis 
 \begin{equation}\label{t}
 \Omega^\top_\lambda(u)=-\Omega_{\lambda^\star}(-u)
 \end{equation}
Combining (\ref{t}) and (\ref{omega1})  we prove  that the shifted determinant $D_\lambda (u)$ is the trace 
 of the composition of two Young symmetrizers, applied to several copies of shifted  matrix $S$. In other words, $D_\lambda(u)$
 is a result of some sort of plethysm of (\ref{Sl}). 
\begin{prop}
Let $\lambda\vdash M$, $\text{dim}\,V_\lambda =(m+1)$. Then 
\begin{equation}\label{pl}
D_\lambda(u) = tr \left(\,\prod_{s=0}^{m} (u-s)\prod_{k=1}^{M} 
\left(\frac{S_{1,sM+k+1}(u+s-m-c_k )}{u-s-c_k }  \right)
 \,(id \otimes \, F_\lambda^{\,\otimes (m+1)} \cdot \, Asym_{m+1})\, \right)
\end{equation}
\end{prop}
\begin{proof}
 Recall from  Section \ref{Sdef} that 
 $$
\alpha\, (\Omega_\lambda(u-m),\, \dots, \, \Omega_\lambda(u))=D_\lambda(u).
$$
Hence, 
\begin{equation*}
\begin{split}
D_\lambda= tr \left(\, \Lambda^{m+1}(\Omega_\lambda(u-m), \,\dots, \, \Omega_\lambda(u))\,\right)\\
= tr \left(\,  ( (\Omega^\top_\lambda(u-m))_{1\, 2}\,\dots \,(\Omega^\top_\lambda(u))_{1 \, m+2})\, Asym_{m+1}\, \right)\\
=  (-1)^{m+1}\,tr\left(\,   (\, (\Omega_{\lambda^{\star}}(-u+m))_{12} \,\dots \, (\Omega_{\lambda^{\star}}(-u))_{1\,m+2})\, Asym_{m+1} \,\right)\\
= tr \left(\,\prod_{s=0}^{m}(u-s) \prod_{k=1}^{M}\frac{1}{(u-s-c_k)}\,( S_\lambda(-u+m))_{1, 2}
\,\dots  \,( S_\lambda(-u))_{1\, m+2}) Asym_{m+1} \,\right)\\
= tr \left(\,\prod_{s=0}^{m} (u-s)\prod_{k=1}^{M} \left(\frac{S_{1,sM+k+1}(u+s-m-c_k)}{u-s-c_k }  \right)
 \,(id \otimes \, F_\lambda^{\,\otimes (m+1)} \cdot \, Asym_{m+1})\, \right)
\end{split}
\end{equation*}
\end{proof}
 We hope that the plethysm relation (\ref{pl}) and the connection of traces of $\Omega_\lambda(u)$ to Capelli elements
  will  allow to prove  Conjecture \ref{conj} in general.

    \end{document}